\input amstex
\input amsppt.sty
\magnification=\magstep1
\hsize=30truecc
\vsize=22.2truecm
\baselineskip=16truept
\NoBlackBoxes
\TagsOnRight \pageno=1 \nologo
\def\Z{\Bbb Z}
\def\N{\Bbb N}

\def\l{\left}
\def\r{\right}
\def\bg{\bigg}
\def\({\bg(}
\def\[{\bg\lfloor}
\def\){\bg)}
\def\]{\bg\rfloor}
\def\t{\text}
\def\f{\frac}

\def\bi{\binom}
\def\eq{\equiv}

\def\ls{\leqslant}
\def\gs{\geqslant}
\def\mo{\roman{mod}}

\def\Proof{\noindent{\it Proof}}

\def\Remark{\medskip\noindent{\it  Remark}}

\hbox {Submitted version, {\tt arXiv:1004.4623v4.}}
\bigskip
\topmatter
\title Products and sums divisible by central binomial coefficients\endtitle
\author Zhi-Wei Sun\endauthor
\leftheadtext{Zhi-Wei Sun} \rightheadtext{Products and sums
divisible by binomial coefficients}
\affil Department of Mathematics, Nanjing University\\
 Nanjing 210093, People's Republic of China
  \\  zwsun\@nju.edu.cn
  \\ {\tt http://math.nju.edu.cn/$\sim$zwsun}
\endaffil
\abstract In this paper we initiate the study of products and sums
divisible by central binomial coefficients. We show that
$$2(2n+1)\bi{2n}n\ \bigg|\ \bi{6n}{3n}\bi{3n}n\ \ \t{for all}\ n=1,2,3,\ldots.$$
Also, for any nonnegative integers $k$ and $n$ we have
$$\bi {2k}k\ \bigg|\ \bi{4n+2k+2}{2n+k+1}\bi{2n+k+1}{2k}\bi{2n-k+1}n$$
and
$$\bi{2k}k\ \bigg|\ (2n+1)\bi{2n}nC_{n+k}\bi{n+k+1}{2k},$$
where $C_m$ denotes the Catalan number
$\f1{m+1}\bi{2m}m=\bi{2m}m-\bi{2m}{m+1}$. On the basis of these results, we obtain certain sums
divisible by central binomial coefficients.
\endabstract
\thanks 2010 {\it Mathematics Subject Classification}.\,Primary 11B65;
Secondary 05A10,  11A07.
\newline\indent Supported by the National Natural Science
Foundation (grant 10871087) and the Overseas Cooperation Fund (grant 10928101) of China.
\endthanks
\endtopmatter
\document

\heading{1. Introduction}\endheading

Central binomial coefficients are given by
$\bi{2n} n$ with  $n\in\N=\{0,1,2,\ldots\}$.
The Catalan numbers
$$C_n=\f1{n+1}\bi{2n}n=\bi{2n}n-\bi{2n}{n+1}\ (n=0,1,2,\ldots)$$
play important roles in combinatorics. (See, e.g., [St].)
There are many sophisticated congruences involving central binomial coefficients and Catalan numbers
(cf. [ST1,ST2] and [S10a,S10b]).

In this paper we investigate a new kind of divisibility problems involving central binomial coefficients.

Our first theorem is as follows.

\proclaim{Theorem 1.1} {\rm (i)} For any positive integer $n$ we have
$$2(2n+1)\bi{2n}n\ \bigg|\ \bi{6n}{3n}\bi{3n}n.\tag1.1$$

{\rm (ii)} Let $k$ and $n$ be nonnegative integers. Then
$$\bi {2k}k\ \bigg|\ \bi{4n+2k+2}{2n+k+1}\bi{2n+k+1}{2k}\bi{2n-k+1}n\tag1.2$$
and
$$\bi{2k}k\ \bigg|\ (2n+1)\bi{2n}nC_{n+k}\bi{n+k+1}{2k}.\tag1.3$$
\endproclaim

In view of (1.1) it is worth introducing the sequence
$$S_n=\f{\bi{6n}{3n}\bi{3n}n}{2(2n+1)\bi{2n}n}\ \ (n=1,2,3,\ldots).$$
Here we list the values of $S_1,\ldots,S_8$:
$$\gather 5,\ 231,\ 14568,\ 1062347,\ 84021990,
\\7012604550,\ 607892634420,\ 54200780036595.
\endgather$$
The author has created this sequence as A176898 at N.J.A Sloane's OEIS (cf. [S10c]).
By Stirling's formula, $S_n\sim 108^n/(8n\sqrt{n\pi})$ as $n\to+\infty$.
 Set $S_0=1/2$.  Using {\tt Mathematica} we find that
  $$\sum_{k=0}^\infty S_kx^k=\f{\sin(\f23\arcsin(6\sqrt{3x}))}{8\sqrt{3x}}\ \ \l(0<x<\f1{108}\r)$$
  and in particular
  $$\sum_{k=0}^\infty\f{S_k}{108^k}=\f{3\sqrt3}8.$$
 {\tt Mathematica} also yields that
 $$\sum_{k=0}^\infty\f{S_k}{(2k+3)108^k}=\f{27\sqrt3}{256}.$$
 It would be interesting to find a combinatorial interpretation or recursion for the sequence $\{S_n\}_{n\gs1}$.

 One can easily show that $S_p\eq15-30p+60p^2\ (\mo\ p^3)$ for any odd prime $p$.
Below we present a conjecture concerning congruence properties of the sequence $\{S_n\}_{n\gs1}$.

\proclaim{Conjecture 1.1} {\rm (i)} Let $n\in\Z^+=\{1,2,3,\ldots\}$. Then
 $S_n$ is odd if and only if $n$ is a power of two. Also, $3S_n\eq0\ (\mo\ 2n+3)$.

 {\rm (ii)} For any prime $p>3$ we have
 $$\sum_{k=1}^{p-1}\f{S_k}{108^k}\eq\cases 0\ (\mo\ p)&\t{if}\ p\eq\pm1\ (\mo\ 12),
 \\-1\ (\mo\ p)&\t{if}\ p\eq\pm5\ (\mo\ 12).\endcases$$
 \endproclaim
 \Remark. Part (i) of Conjecture 1.1 might be
 shown by our method for proving Theorem 1.1(i), but we are
 not interested in writing the details.

 \medskip

Our following conjecture is concerned with a companion sequence of $\{S_n\}_{n\gs0}$.

\proclaim{Conjecture 1.2} There are positive integers $T_1,T_2,T_3,\ldots$ such that
$$\sum_{k=0}^\infty S_kx^{2k+1}+\f1{24}-\sum_{k=1}^\infty T_kx^{2k}
=\f{\cos(\f23\arccos(6\sqrt 3x))}{12}$$
for all real $x$ with $|x|\ls 1/(6\sqrt3)$.
Also, $T_p\eq-2\ (\mo\ p)$ for any prime $p$.
\endproclaim
Here we list the values of $T_1,\ldots,T_8$:
$$\gather 1,\ 32,\ 1792,\ 122880,\ 9371648,
\\763363328,\ 65028489216,\ 5722507051008.
\endgather$$

In 1914 Ramanujan [R] obtained that
$$\sum_{k=0}^\infty\f{4k+1}{(-64)^k}\bi{2k}k^3=\f2{\pi}$$
and
$$\sum_{k=0}^\infty(20k+3)\f{\bi{2k}k^2\bi{4k}{2k}}{(-2^{10})^k}=\f 8{\pi}.$$
(See also [BB], [BBC] and B. C. Berndt [Be] for such series.)
Actually the first identity was originally proved by G. Bauer in
1859.
 Both identities can be
proved via the WZ (Wilf-Zeilberger) method (see M. Petkov\v sek, H.
S. Wilf and D. Zeilberger [PWZ], and Zeilberger [Z] for this
method),  for example,
 Guillera [G] used the WZ method to prove the second identity.
van Hammer [vH] conjectured that the first identity has a $p$-adic analogue.
This conjecture was first proved by E. Mortenson [M], and recently re-proved in [Zu] via the WZ method.

On the basis of Theorem 1.1, we deduce the following result
which was conjectured by the author in [S10b].

\proclaim{Theorem 1.2} For any positive integer $n$ we have
$$4(2n+1)\bi{2n}n\ \bigg|\ \sum_{k=0}^{n}(4k+1)\bi{2k}k^3(-64)^{n-k}\tag1.4$$
and
$$4(2n+1)\bi{2n}n\ \bigg|\ \sum_{k=0}^{n}(20k+3)\bi{2k}k^2\bi{4k}{2k}(-2^{10})^{n-k}.\tag1.5$$
\endproclaim
\Remark. In 1998 N. J. Calkin [C] proved that $\bi{2n}n\mid\sum_{k=-n}^n(-1)^k\bi{2n}{n+k}^m$
for any $m,n\in\Z^+$. See also V.J.W. Guo, F. Jouhet and J. Zeng [GJZ], and H.Q. Cao and H. Pan [CP]
for further extensions of Calkin's result.

\medskip

Now we raise two more conjectures.

\proclaim{Conjecture 1.3} {\rm (i)} For any $n\in\Z^+$ we have
$$a_n:=\f1{8n^2\bi{2n}n^2}\sum_{k=0}^{n-1}(205k^2+160k+32)(-1)^{n-1-k}\bi{2k}k^5\in\Z^+.$$

{\rm (ii)} Let $p$ be an odd prime. If $p\not=3$ then
$$\sum_{k=0}^{(p-1)/2}(205k^2+160k+32)(-1)^k\bi{2k}k^5\eq 32p^2+\f{896}3p^5B_{p-3}\ (\mo\ p^6),$$
where $B_0,B_1,B_2,\ldots$ are Bernoulli numbers.
If $p\not=5$ then
$$\sum_{k=0}^{p-1}(205k^2+160k+32)(-1)^k\bi{2k}k^5\eq 32p^2+64p^3H_{p-1}\ (\mo\ p^7),$$
where $H_{p-1}=\sum_{k=1}^{p-1}1/k$.
\endproclaim
\Remark. Note that $a_1=1$ and
$$4(2n+1)^2a_{n+1}+n^2a_n=(205n^2+160n+32)\bi{2n-1}n^3\ \ \ \t{for}\ n=1,2,\ldots.$$
The author created the sequence $\{a_n\}_{n>0}$ at OEIS as A176285 (cf. [S10c]).
In 1997 T. Amdeberhan and D. Zeilberger [AZ] used the WZ method to obtain
$$\sum_{k=1}^\infty\f{(-1)^k(205k^2-160k+32)}{k^5\bi{2k}k^5}=-2\zeta(3).$$

\proclaim{Conjecture 1.4} {\rm (i)} For any odd prime $p$, we have
$$\sum_{k=0}^{p-1}\f{28k^2+18k+3}{(-64)^k}\bi{2k}k^4\bi{3k}k\eq 3p^2-\f 72p^5B_{p-3}\ (\mo\ p^6),$$
and
$$\sum_{k=0}^{(p-1)/2}\f{28k^2+18k+3}{(-64)^k}\bi{2k}k^4\bi{3k}k\eq 3p^2+6\l(\f{-1}p\r)p^4E_{p-3}\ (\mo\ p^5),$$
where $E_0,E_1,E_2,\ldots$ are Euler numbers.

{\rm (ii)} For any integer $n>1$, we have
$$\sum_{k=0}^{n-1}(28k^2+18k+3)\bi{2k}k^4\bi{3k}k(-64)^{n-1-k}\eq0\ \ \(\mo\ (2n+1)n^2\bi{2n}n^2\).$$
Also,
$$\sum_{k=1}^\infty\f{(28k^2-18k+3)(-64)^k}{k^5\bi{2k}k^4\bi{3k}k}=-14\zeta(3).$$
\endproclaim
\Remark. The conjectured series for $\zeta(3)=\sum_{n=1}^\infty1/n^3$ was first announced by the author
in a message to Number Theory Mailing List (cf. [S10d]) on April 4, 2010.

\medskip

For more conjectures similar to Conjectures 1.3 and 1.4 the reader may consult [S09] and [S10c].

In the next section we will establish three auxiliary inequalities involving the floor function.
Sections 3 and 4 are devoted to the proofs of Theorem 1.1 and Theorem 1.2 respectively.

\heading{2. Three auxiliary inequalities}\endheading

In this section, for a rational number $x$ we let $\{x\}=x-\lfloor
x\rfloor$ be the fractional part of $x$, and set $\{x\}_m=m\{x/m\}$
for any $m\in\Z^+$.

\proclaim{Theorem 2.1}  Let $m>1$ be an integer. Then for any $n\in\Z$ we have
$$\l\lfloor\f nm\r\rfloor+\l\lfloor\f{6n}m\r\rfloor\gs\l\lfloor\f{2n}m\r\rfloor+\l\lfloor\f{2n+1}m\r\rfloor
+\l\lfloor\f{3n}m\r\rfloor.\tag2.1$$
\endproclaim
\Proof. Let $A_m(n)$ denote the left-hand side of (2.1) minus the
right-hand side. Then
$$A_m(n)=\l\{\f{2n}m\r\}+\l\{\f{2n+1}m\r\}+\l\{\f{3n}m\r\}-\f1m-\l\{\f nm\r\}-\l\{\f{6n}m\r\},$$
which only depends on $n$ modulo $m$. So, without any loss of generality we may simply assume that $n\in\{0,\ldots,m-1\}$.
Hence $A_m(n)\gs0$ if and only if
$$\l\{\f{2n}m\r\}+\l\{\f{2n+1}m\r\}+\l\{\f{3n}m\r\}\gs \f{n+1}m.\tag2.2$$
(Note that $2n+(2n+1)+3n-(n+1)=6n$.)

(2.1) is obvious when $n=0$. If $1\ls n<m/2$, then $\{2n/m\}=2n/m\gs(n+1)/m$ and hence (2.2) holds.
In the case $n\gs m/2$, (2.2) can be simplified as
$$\f{3n}m+\l\{\f{3n}m\r\}\gs2,$$
which holds since $3n\gs m+m/2$.

By the above we have proved (2.1).
\qed

\proclaim{Theorem 2.2} Let $m\in\Z^+$ and $k,n\in\Z$. Then we have
$$\l\lfloor\f{4n+2k+2}m\r\rfloor-\l\lfloor\f{2n+k+1}m\r\rfloor+2\l\lfloor \f km\r\rfloor-2\l\lfloor\f{2k}m\r\rfloor
\gs\l\lfloor\f nm\r\rfloor+\l\lfloor\f{n-k+1}m\r\rfloor,\tag2.3$$
unless $2\mid m$ and $k\eq n+1\eq m/2\ (\mo\ m)$ in which case the
right-hand side of the inequality equals the left-hand side plus
one.
\endproclaim
\Proof. Since
$$(4n+2k+2)-(2n+k+1)+2k-2(2k)=n+(n-k+1),$$
(2.3) has the following equivalent form:
$$\l\{\f{4n+2k+2}m\r\}-\l\{\f{2n+k+1}m\r\}+2\l\{\f
km\r\}-2\l\{\f{2k}m\r\}\ls \l\{\f nm\r\}+\l\{\f{n-k+1}m\r\}.\tag2.4$$
Note that this only depends on $k$ and $n$ modulo $m$. So, without any loss of generality,
we may simply assume that $k,n\in\{0,\ldots,m-1\}$.

{\it Case} 1. $k<m/2$ and $\{2n+k+1\}_m<m/2$.

In this case, (2.4) can be simplified as
$$\f{n+2k}m+\l\{\f{n-k+1}m\r\}\gs\l\{\f{2n+k+1}m\r\},$$
which is true since the left-hand side is nonnegative and $(n+2k)+(n-k+1)\eq 2n+k+1\ (\mo\ m)$.

{\it Case} 2. $k<m/2$ and $\{2n+k+1\}_m\gs m/2$.

In this case, (2.4) can be simplified as
$$\f{n+2k}m+\l\{\f{n-k+1}m\r\}\gs\l\{\f{2n+k+1}m\r\}-1,$$
which holds trivially since the right-hand side is negative.

{\it Case} 3. $k\gs m/2$ and $\{2n+k+1\}_m<m/2$.

In this case, (2.4) can be simplified as
$$\f{n+2k}m+\l\{\f{n-k+1}m\r\}\gs 2+\l\{\f{2n+k+1}m\r\}.$$
Since $(n+2k)+(n-k+1)=2n+k+1$, this is equivalent to
$$n+2k+\{n-k+1\}_m\gs 2m.$$

If $k>n+1$, then
$$n+2k+\{n-k+1\}_m=n+2k+(n-k+1+m)=2n+k+1+m\gs 2m$$
since $2n+k+1>k\gs m/2$ and $\{2n+k+1\}_m<m/2$.

Now assume that $k\ls n+1$. Clearly
$$n+2k+\{n-k+1\}_m=n+2k+(n-k+1)=2n+k+1\gs 3k-1.$$
If $k>m/2$ then $3k-1\gs 3(m+1)/2-1>3m/2$.
If $k\ls n$ then $2n+k+1>3k\gs 3m/2$.
So, except the case $k=n+1=m/2$ we have
$$n+2k+\{n-k+1\}_m=2n+k+1\gs 3m/2$$ and
hence $n+2k+\{n-k+1\}_m=2n+k+1\gs 2m$ since $\{2n+k+1\}_m<m/2$.

When $k=n+1=m/2$, the left-hand side of (2.4) minus the right-hand side equals
$$\f{m-2}m-\f{m/2-1}m+2\f{m/2}m-\f{m/2-1}m=1.$$

{\it Case} 4. $k\gs m/2$ and $\{2n+k+1\}_m\gs m/2$.

In this case, clearly $m\not=1$, and (2.4) can be simplified as
$$\f{n+2k}m+\l\{\f{n-k+1}m\r\}\gs 1+\l\{\f{2n+k+1}m\r\}$$
which is equivalent to
$$n+2k+\{n-k+1\}_m\gs m.$$
If $k\ls n+1$, then
$$n+2k+\{n-k+1\}_m=n+2k+(n+1-k)=2n+k+1\gs 3k-1\gs \f {3m}2-1\gs m.$$
If $k>n+1$, then
$$n+2k+\{n-k+1\}_m=n+2k+(n+1-k)+m=2n+k+1+m>m.$$

In view of the above, we have completed the proof of Theorem 2.2.

\proclaim{Theorem 2.3} Let $m\in\Z^+$ and $k,n\in\Z$. Then we have
$$\aligned&\l\lfloor\f{2n+2k}m\r\rfloor-\l\lfloor\f{n+k}m\r\rfloor+2\l\lfloor \f km\r\rfloor-2\l\lfloor\f{2k}m\r\rfloor
\\&\ \ \gs2\l\lfloor\f
nm\r\rfloor-\l\lfloor\f{2n+1}m\r\rfloor+\l\lfloor\f{n-k+1}m\r\rfloor,
\endaligned\tag2.5$$
unless $2\mid m$ and $k\eq n+1\eq m/2\ (\mo\ m)$ in which case the
right-hand side of the inequality equals the left-hand side plus
one.
\endproclaim
\Proof. Since
$$2n+2k-(n+k)+2k-2(2k)=2n-(2n+1)+(n-k+1),$$
(2.5) is equivalent to the following inequality:
$$\aligned&\l\{\f{2n+2k}m\r\}-\l\{\f{n+k}m\r\}+2\l\{\f km\r\}-2\l\{\f{2k}m\r\}
\\ &\ \ \ls2\l\{\f nm\r\}-\l\{\f{2n+1}m\r\}+\l\{\f{n-k+1}m\r\}.
\endaligned\tag2.6$$
As (2.6) only depends on $k$ and $n$ modulo $m$, without loss of generality we simply assume that $k,n\{0,\ldots,m-1\}$.

{\it Case} 1. $k<m/2$ and $\{n+k\}_m<m/2$.

In this case, (2.6) can be simplified as
$$\f{2n+2k}m+\l\{\f{n-k+1}m\r\}\gs\l\{\f{2n+1}m\r\}+\l\{\f{n+k}m\r\}$$
which holds since
$$\f{2n+2k}m-\l\{\f{n+k}m\r\}+\l\{\f{n-k+1}m\r\}\gs0$$
and $2n+2k-(n+k)+(n-k+1)=2n+1.$

{\it Case} 2. $k<m/2$ and $\{n+k\}_m\gs m/2$.

In this case, (2.6) can be simplified as
$$\f{2n+2k}m+\l\{\f{n-k+1}m\r\}\gs\l\{\f{2n+1}m\r\}+\l\{\f{n+k}m\r\}-1$$
which holds since
$$\f{2n+2k}m\gs\f{n+k}m\gs\l\{\f{n+k}m\r\}
\ \t{and}\ \l\{\f{n-k+1}m\r\}\gs0>\l\{\f{2n+1}m\r\}-1.$$

{\it Case} 3. $k\gs m/2$ and $\{n+k\}_m<m/2$.

In this case, we must have $n+k\gs m$ and hence $\{n+k\}_m=n+k-m$. Thus (2.6) can be simplified as
$$\f{n+k-m}m+\l\{\f{n-k+1}m\r\}\gs\l\{\f{2n+1}m\r\}$$
which holds trivially since $n+k-m+(n-k+1)\eq 2n+1\ (\mo\ m)$.

{\it Case} 4. $k\gs m/2$ and $\{n+k\}_m\gs m/2$.

In this case, (2.6) can be simplified as
$$\f{2n+2k}m-\l\{\f{n+k}m\r\}+\l\{\f{n-k+1}m\r\}\gs 1+\l\{\f{2n+1}m\r\}$$
which is equivalent to
$$\f{2(n+k)}m-\l\{\f{n+k}m\r\}+\l\{\f{n-k+1}m\r\}\gs1\tag2.7$$
since $2n+2k-(n+k)+(n-k+1)=2n+1$.

Clearly (2.7) holds if $n+k\gs m$.
If $n+k<m$ and $k>n+1$, then the left-hand side of the inequality (2.7) is
$$\f{n+k}m+\f{n+1-k}m+1=\f{2n+1}m+1>1.$$

Now assume that $n+k<m$ and $k\ls n+1$. Then (2.7) is equivalent to $2n+1\gs m$.
If $k\ls n$ then $2n+1>2k\gs m$. If $k=n+1\not= m/2$, then
$k=n+1\gs(m+1)/2$ and hence $2n+1=2(n+1)-1\gs m$.

When $k=n+1=m/2$, the left-hand side of (2.6) minus the right-hand side equals
$$\f{m-2}m-\f{m-1}m+2\f{m/2}m-2\f{m/2-1}m+\f{m-1}m=1.$$

Combining the discussion of the four cases we obtain the desired result. \qed

\heading{3. Proof of Theorem 1.1}\endheading

For a prime $p$,  the $p$-adic evaluation of an integer $m$ is given by
$$\nu_p(m)=\sup\{a\in\N:\ p^a\mid m\}.$$
For a rational number $x=m/n$ with $m\in\Z$ and $n\in\Z^+$, we set $\nu_p(x)=\nu_p(m)-\nu_p(n)$ for any prime $p$.
Note that a rational number $x$ is an integer if and only if $\nu_p(x)\gs0$ for all primes $p$.

\medskip\noindent
{\it Proof of Theorem 1.1}. (i) Fix $n\in\Z^+$, and define $A_m(n)$ for $m>1$ as in the proof of Theorem 2.1.
Observe that
$$Q:=\f{\bi{6n}{3n}\bi{3n}n}{(2n+1)\bi{2n}n}=\f{n!(6n)!}{(2n)!(2n+1)!(3n)!}.$$
So, for any prime $p$ we have
$$\nu_p(Q)=\sum_{i=1}^\infty A_{p^i}(n)\gs0$$
by Theorem 2.1. Therefore $Q$ is an integer.

Choose $j\in\Z^+$ such that $2^{j-1}\ls n<2^j$.
As $2n+1\ls 2(2^j-1)+1<2^{j+1}$, we have
$$\align &\l\lfloor\f n{2^{j+1}}\r\rfloor+\l\lfloor\f{6n}{2^{j+1}}\r\rfloor
-\l\lfloor\f{2n}{2^{j-1}}\r\rfloor-\l\lfloor\f{2n+1}{2^{j-1}}\r\rfloor
-\l\lfloor\f{3n}{2^{j-1}}\r\rfloor
\\=&\l\lfloor\f{3n}{2^j}\r\rfloor-\l\lfloor\f{3n}{2^{j+1}}\r\rfloor
=\l\lfloor\f{3n+2^j}{2^{j+1}}\r\rfloor\gs\l\lfloor\f{2n+2^j}{2^{j+1}}\r\rfloor\gs1.
\endalign$$
Therefore
$$\nu_2(Q)=\sum_{i=1}^\infty A_{2^i}(n)\gs A_{2^{j+1}}(n)\gs1.$$
and hence $Q$ is even. This proves (1.1). \qed

(ii) (1.2) and (1.3) are obvious in the case $k=0$. If $k>n+1$, then
$$\bi{2n+k+1}{2k}=\bi{n+k+1}{2k}=0$$
and hence (1.2) and (1.3) hold trivially.
Below we assume that $1\ls k\ls n+1$.

 Recall that for any nonnegative integer $m$ and prime $p$ we have
 $$\nu_p(m!)=\sum_{i=1}^\infty\l\lfloor\f m{p^i}\r\rfloor.$$
 Since
 $$\f{\bi{4n+2k+2}{2n+k+1}\bi{2n+k+1}{2k}\bi{2n+k+1}n}{\bi{2k}k}
 =\f{(4n+2k+2)!(k!)^2}{(2n+k+1)!((2k)!)^2n!(n-k+1)!}$$
 and
 $$\f{(2n+1)\bi{2n}nC_{n+k}\bi{n+k+1}{2k}}{\bi{2k}k}=\f{(2n+1)!(2n+2k)!(k!)^2}{(n!)^2(n+k)!((2k)!)^2(n-k+1)!},$$
 it suffices to show that for any prime $p$ we have
 $$\sum_{i=1}^\infty C_{p^i}(n,k)\gs0\ \ \t{and}\ \ \sum_{i=1}^\infty D_{p^i}(n,k)\gs0,$$
 where
 $$\align C_m(n,k)=&\l\lfloor\f{4n+2k+2}m\r\rfloor-\l\lfloor\f{2n+k+1}m\r\rfloor
 +2\l\lfloor\f km\r\rfloor-2\l\lfloor\f{2k}m\r\rfloor
 \\&-\l\lfloor\f nm\r\rfloor-\l\lfloor\f{n-k+1}m\r\rfloor
 \endalign$$
 and
 $$\align D_m(n,k)=&\l\lfloor\f{2n+2k}m\r\rfloor-\l\lfloor\f{n+k}m\r\rfloor+2\l\lfloor\f km\r\rfloor-2\l\lfloor\f{2k}m\r\rfloor
 \\&-2\l\lfloor\f nm\r\rfloor+\l\lfloor\f{2n+1}m\r\rfloor-\l\lfloor\f{n-k+1}m\r\rfloor.
\endalign$$

(a) By Theorem 2.2, $C_{p^i}(n,k)\gs0$ unless $p=2$ and $k\eq n+1\eq 2^{i-1}\ (\mo\ 2^i)$
in which case $C_{2^i}(n,k)=-1$. Suppose that $k\eq n+1\eq 2^{i-1}\ (\mo\ 2^i)$, $k=2^{i-1}k_0$ and $n+1=2^{i-1}n_0$, where
$1\ls k_0\ls n_0$ and $k_0$ and $n_0$ are odd. If $i\gs2$, then
$$C_{2^{i-1}}(n,k)=4n_0+2k_0-1-(2n_0+k_0-1)+2k_0-4k_0-(n_0-1)-(n_0-k_0)=1$$
and hence $C_{2^{i-1}}(n,k)+C_{2^i}(n,k)=1+(-1)=0$.
So it remains to consider the case $k\eq n+1\eq1\ (\mo\ 2)$.

Assume that $k$ is odd and $n$ is even. Write $k+1=2^jk_1$ and $n=2n_1$ with $k_1,n_1\in\Z^+$ and $2\nmid k_1$.
Then it is easy to check that
$$\align C_{2^{j+1}}(n,k)=&\l\lfloor\f{4n_1}{2^j}\r\rfloor+k_1-\l\lfloor\f{2n_1-2^{j-1}
+2^{j-1}(k_1-1)}{2^j}\r\rfloor
\\&+2\l\lfloor\f{k_1}2\r\rfloor-2\l\lfloor\f{2^jk_1-1}{2^j}\r\rfloor
-\l\lfloor\f{n_1}{2^j}\r\rfloor-\l\lfloor\f{n_1+1-2^{j-1}k_1}{2^j}\r\rfloor
\\=&\l\lfloor\f{4n_1}{2^j}\r\rfloor+k_1-\l\lfloor\f{2n_1-2^{j-1}}{2^j}\r\rfloor-\f{k_1+1}2+k_1-1-2(k_1-1)
\\&-\l\lfloor\f{n_1}{2^j}\r\rfloor-\l\lfloor\f{n_1+1+2^{j-1}}{2^j}\r\rfloor+\f{k_1+1}2
\\=&1+\l\lfloor \f{n_1+(n_1+1+2^{j-1})+(2n_1-2^{j-1})}{2^j}\r\rfloor
\\&-\l\lfloor\f {n_1}{2^j}\r\rfloor-\l\lfloor\f{n_1+1+2^{j-1}}{2^j}\r\rfloor-\l\lfloor\f{2n_1-2^{j-1}}{2^j}\r\rfloor
\\\gs&1
\endalign$$
and hence $C_2(n,k)+C_{2^{j+1}}(n,k)\gs0$.

By the above, we do have $\sum_{i=1}^\infty C_{p^i}(n,k)\gs0$ for any prime $p$. So (1.2) holds.

(b) By Theorem 2.2, $D_{p^i}(n,k)\gs0$ unless $p=2$ and $k\eq n+1\eq 2^{i-1}\ (\mo\ 2^i)$
in which case $D_{2^i}(n,k)=-1$. So, to prove (1.2) it suffices to find a positive integer $j$
such that $D_{2^j}(n,k)\gs1$.

Clearly there is a unique positive integer $j$ such that $2^{j-1}\ls n+k<2^j$. Note that $k\ls(n+k)/2<2^{j-1}$
and $$D_{2^j}(n,k)=1+\l\lfloor\f{2n+1}{2^j}\r\rfloor\gs1.$$
This concludes the proof of (1.3).

  The proof of Theorem 1.1 is now complete. \qed

\heading{4. Proof of Theorem 1.2}\endheading

\medskip
\noindent{\it Proof of Theorem 1.2}. (i) We first prove (1.4). For $k,n\in\N$ define
$$F(n,k)=\f{(-1)^{n+k}(4n+1)}{4^{3n-k}}\bi{2n}n^2\f{\bi{2n+2k}{n+k}\bi{n+k}{2k}}{\bi{2k}k}$$
and
$$G(n,k)=\f{(-1)^{n+k}(2n-1)^2\bi{2n-2}{n-1}^2}{2(n-k)4^{3(n-1)-k}}\bi{2(n-1+k)}{n-1+k}\f{\bi{n-1+k}{2k}}{\bi{2k}k}.$$
Clearly $F(n,k)=G(n,k)=0$ if $n<k$. By [Zu],
$$F(n,k-1)-F(n,k)=G(n+1,k)-G(n,k)$$
for all $k\in\Z^+$ and $n\in\N$.

Fix a positive integer $N$. Then
$$\align \sum_{n=0}^NF(n,0)-F(N,N)=&\sum_{n=0}^NF(n,0)-\sum_{n=0}^NF(n,N)
\\=&\sum_{k=1}^N\(\sum_{n=0}^NF(n,k-1)-\sum_{n=0}^NF(n,k)\)
\\=&\sum_{k=1}^N\sum_{n=0}^N(G(n+1,k)-G(n,k))
=\sum_{k=1}^NG(N+1,k).
\endalign$$
Note that
$$\sum_{n=0}^NF(n,0)=\sum_{n=0}^N\f{4n+1}{(-64)^n}\bi{2n}n^3$$
and
$$F(N,N)=\f{4N+1}{4^{2N}}\bi{2N}N\bi{4N}{2N}=\f{(4N+1)(2N+1)}{4^{2N}}\bi{2N}NC_{2N}.$$
Also,
$$\align\sum_{k=1}^NG(N+1,k)=&\f{(2N+1)^2}2\sum_{k=1}^N\f{(-1)^{N+k+1}}{4^{3N-k}}\bi{2N}N^2C_{N+k}\f{\bi{N+k+1}{2k}}{\bi{2k}k}
\\=&\f{2(2N+1)\bi{2N}N}{(-64)^N}\sum_{k=1}^N(-4)^{k-1}\f{(2N+1)\bi{2N}NC_{N+k}\bi{N+k+1}{2k}}{\bi{2k}k}.
\endalign$$
and
$$\align \f{\bi{2N}NC_{N+1}\bi{N+2}2}{\bi 21}=&\bi{2N-1}{N-1}\bi{2N+2}{N+1}\f{N+1}2
\\=&\bi{2N-1}{N-1}\bi{2N+1}{N+1}(N+1)
\\=&\bi{2N-1}{N-1}(2N+1)\bi{2N}N
\\=&2(2N+1)\bi{2N-1}{N-1}^2\eq0\ (\mo\ 2).
\endalign$$
So, with the help of (1.3) we see that $\sum_{n=0}^N(4n+1)\bi{2n}n^3(-64)^{N-n}$
is divisible by $4(2N+1)\bi{2N}N$.

(ii) Now we turn to the proof of (1.5).

For $n,k\in\N$, define
$$F(n,k):=\f{(-1)^{n+k}(20n-2k+3)}{4^{5n-k}}\cdot\f{\bi{2n}n\bi{4n+2k}{2n+k}\bi{2n+k}{2k}\bi{2n-k}n}{\bi{2k}k}.
$$
and
$$G(n,k)
:=\f{(-1)^{n+k}}{4^{5n-4-k}}\cdot\f{n\bi{2n}n\bi{4n+2k-2}{2n+k-1}\bi{2n+k-1}{2k}\bi{2n-k-1}{n-1}}{\bi{2k}k}.$$
Clearly $F(n,k)=G(n,k)=0$ if $n<k$. By [Zu],
$$F(n,k-1)-F(n,k)=G(n+1,k)-G(n,k)$$
for all $k\in\Z^+$ and $n\in\N$.

Fix a positive integer $N$. As in part (i) we have
$$\sum_{n=0}^NF(n,0)-F(N,N)=\sum_{k=1}^NG(N+1,k).$$
Observe that
$$\sum_{n=0}^NF(n,0)=\sum_{n=0}^N\f{20n+3}{(-2^{10})^n}\bi{2n}n^2\bi{4n}{2n}$$
and
$$F(N,N)=\f{18N+3}{2^{8N}}\bi{6N}{3N}\bi{3N}N.$$
Also,
$$\sum_{k=1}^NG(N+1,k)
=\f{2(2N+1)\bi{2N}N}{(-2^{10})^N}\sum_{k=1}^N(-4)^{k-1}\f{\bi{4N+2k+2}{2N+k+1}\bi{2N+k+1}{2k}\bi{2N-k+1}N}{\bi{2k}k}.$$
Note that
$$\f{\bi{4N+4}{2n+2}\bi{2N+2}2\bi{2N}N}{\bi 21}=2\bi{4N+3}{2N+1}\bi{2N+2}2\bi{2N-1}{N-1}\eq0\ (\mo\ 2).$$
Applying (1.2) we see that $(-2^{10})^N\sum_{k=1}^NG(N+1,k)$ is a multiple of
$4(2N+1)\bi{2N}N$. By (1.1),
$$(-2^{10})^N\f{18N+3}{2^{8N}}\bi{6N}{3N}\bi{3N}N$$
is divisible by $8(2N+1)\bi{2N}N$.
Therefore
$$\sum_{n=0}^N(20n+3)\bi{2n}n^2\bi{4n}{2n}(-2^{10})^{N-n}$$
is a multiple of $4(2N+1)\bi{2N}N$.
\medskip

Combining the above, we have completed the proof of Theorem 1.2. \qed

\enddocument

\bye